\documentclass[a4paper,11pt]{amsart}

\usepackage{t1enc}
\usepackage[latin1]{inputenc}
\usepackage{graphics}
\usepackage{graphicx}
\usepackage{amsfonts}
\usepackage{eufrak}
\usepackage{xypic}
\usepackage{amssymb}
\usepackage{amsmath}
\usepackage{amsthm}
\usepackage{mathrsfs}

\theoremstyle{plain}
\newtheorem{thm}{Theorem}[section]
\newtheorem{lemma}[thm]{Lemma}
\newtheorem{prop}[thm]{Proposition}
\newtheorem{cor}[thm]{Corollary}

\theoremstyle{definition}
\newtheorem{df}[thm]{Definition}
\newtheorem{rem}[thm]{Remark}

\def\Qp{\mathbb{Q}_p}
\def\Zp{\mathbb{Z}_p}
\def\Hh{\mathcal{H}}
\def\CC{\mathbb{C}}

\def\sg{\mathscr{G}}
\def\sh{\mathscr{H}}
\def\sk{\mathscr{K}}
\def\sn{\mathscr{N}}
\def\sq{\mathscr{Q}}
\providecommand{\gm}{\Gamma}
\providecommand{\h}{\mathcal{H}}

\numberwithin{equation}{section}

\begin{document}

\title[Positive definite $*$-spherical functions]{Positive definite $*$-spherical functions,  property (T), and $C^*$-completions of Gelfand pairs}

\author{Nadia S.~Larsen}
\address{Department of Mathematics \\ University of Oslo \\PO BOX 1053 Blindern\\ Norway}
\email{nadiasl@math.uio.no}

\author{Rui ~Palma}
\email{ruip@math.uio.no}

\thanks{Research partially supported by the Research Council of Norway and the NordForsk research network ``Operator Algebra and Dynamics'' (grant 11580).}

\begin{abstract} The study of existence of a universal $C^*$-completion of the $^*$-algebra canonically associated to a Hecke pair was initiated by Hall, who proved that the Hecke algebra associated to $(\operatorname{SL}_2(\Qp), \operatorname{SL}_2(\Zp))$ does not admit a universal $C^*$-completion. Kaliszewski, Landstad and Quigg studied the problem by placing it in the framework of Fell-Rieffel  equivalence, and highlighted the role of other $C^*$-completions.  In the case of  the pair $(\operatorname{SL}_n(\Qp), \operatorname{SL}_n(\Zp))$ for $n\geq 3$ we show, invoking property (T) of $\operatorname{SL}_n(\Qp)$, that the $C^*$-completion of the $L^1$-Banach algebra and the corner of $C^*(\operatorname{SL}_n(\Qp))$ determined by the subgroup are distinct. In fact, we prove  a more general result valid for a simple algebraic group  of rank at least $2$ over a $\mathfrak{p}$-adic field with a good choice of a maximal compact open subgroup.
\end{abstract}

\date{11 July, 2014. Revised 12 November, 2014 and 26 May, 2015.}

\maketitle

\section{Introduction}\label{sec: intro}
The work of Bost and Connes on a $C^*$-dynamical system with deep connections to class field theory brought to attention  $C^*$-algebras associated to Hecke pairs, \cite{BC}. A (discrete) Hecke pair $(G, \Gamma)$ consists of a group $G$ with a subgroup $\Gamma$ such that every double coset $\Gamma g\Gamma$ contains finitely many left cosets, for every $g$ in $G$. The associated Hecke algebra $\mathcal{H}(G, \Gamma)$ is a convolution $^*$-algebra of complex-valued functions on the space of double cosets, see for example \cite{kri}. Hall initiated the study of $C^*$-completions of $\mathcal{H}(G, \Gamma)$ in connection to asking whether there is a category equivalence between the nondegenerate $^*$-representations of $\mathcal{H}(G, \Gamma)$ and the unitary representations of $G$ generated by their $\Gamma$-fixed vectors, \cite{hall}. The question was motivated by the well-known fact that  unitary representations of a group are in bijective correspondence with the nondegenerate $^*$-representations of the group algebra.

 While proposing a condition that would yield an affirmative answer to the question,  Hall showed at the same time that the equivalence could not hold in all generality. Indeed, she showed that  $\mathcal{H}(\operatorname{SL}_2(\Qp), \operatorname{SL}_2(\Zp))$ does not admit an enveloping $C^*$-algebra, \cite{hall}. The proof was based on a careful description, via the Satake isomorphism, of this algebra as a polynomial algebra in one variable. This showed that there are elements of $\mathcal{H}(\operatorname{SL}_2(\Qp), \operatorname{SL}_2(\Zp))$ which have arbitrarily large norms with respect to $^*$-representations.

Shortly afterwards there were two related developments. In one of them Tzanev, motivated by amenability in connection with Hecke pairs, and drawing on work of Schlichting \cite{schl}, noticed that to every discrete Hecke pair $(G, \Gamma)$ one can associate in a canonical way a topological pair $(\overline{G}, \overline{\Gamma})$ having isomorphic Hecke  algebra, \cite{tzanev}.  In another development, Kaliszewski, Landstad and Quigg placed the study of Hall's correspondence in the framework of Fell-Rieffel equivalence, \cite{kal-land-qui}. Both approaches investigated $C^*$-completions of $\mathcal{H}(G, \Gamma)$ by looking at the corner $p_0C_c(\overline{G})p_0$ arising from the self-adjoint projection $p_0$ given by the characteristic function of the compact open  subgroup $\overline{\Gamma}$ of $\overline{G}$.

One question raised in \cite{kal-land-qui} was when a certain canonical surjection from $C^*(p_0L^1(\overline{G})p_0)$ onto $p_0C^*(\overline{G})p_0$ could fail to  be injective. It is stated in \cite[page 677]{kal-land-qui} that the injectivity fails for the pair $(\operatorname{PSL}_3(\Qp), \operatorname{PSL}_3(\Zp))$, according to a personal communication by Tzanev. However, to our knowledge, no proof of this claim has been published. The result showing that injectivity of the canonical map above fails for the pair $(\operatorname{SL}_2(\Qp), \operatorname{SL}_2(\Zp))$ is due to the second named author, \cite{palma1}.

 Our contribution here is to show that the canonical surjection from $C^*(p_0L^1(\operatorname{SL}_n(\Qp))p_0)$ onto $p_0C^*(\operatorname{SL}_n(\Qp))p_0$ is not injective for all $n\geq 3$. Indeed, we prove the result in much greater generality.  The proof differs from the case of $n=2$ in \cite{palma1} in that it appeals to property (T). After the second named author's talk on the case $n=2$ at a NordForsk conference on the Far{\o}e Islands in 2012, a discussion with M. Landstad revealed that the obstruction to injectivity in the case $n=3$ would be property (T).

 Thus our initial investigations concentrated on property (T) for a Hecke pair $(G, \Gamma)$. This notion was already introduced in  \cite[page 24]{tzanev-thesis}, where it was claimed that one could easily show that property (T) of a Hecke pair $(G, \Gamma)$ was equivalent to property (T) of the locally compact group $\overline{G}$, and that all other characterisations of property (T) for groups could be translated to Hecke pairs. However, no details were given for any of these two claims.  Here we show that the first claim is valid, though the arguments require a careful analysis; this material forms the content of appendix A. We point out that this version of property (T) for Hecke pairs, termed co-rigidity, is also discussed in \cite[\S 6]{AD}. Regarding Tzanev's second claim, we found that one characterisation is more subtle. It is known that for a locally compact group $G$ with property (T), the trivial representation of $G$ is isolated in $\widehat{G}$ with its natural Fell topology, and therefore the trivial representation of $C^*(G)$ is isolated in the hull-kernel topology. However, in taking a compact open subgroup $U$ and forming the Hecke pair $(G, U)$, we found that the trivial representation of $C^*(L^1(G, U))$ need not be isolated in the hull-kernel topology on $(C^*(L^1(G, U)))^\wedge$. In fact, this is precisely what happens when $(G, U)$ is a Gelfand pair formed of an algebraic group over a $\mathfrak{p}$-adic field together with a good choice of a compact open subgroup, see Theorem~\ref{SLn trivial representation not isolated}.

As a consequence of this distinction, the trivial representation is isolated in the corner $p_0C^*(G)p_0$ with its hull-kernel topology, see Proposition~\ref{prop:isolation-passes-to-corner}, but not in $(C^*(L^1(G, U)))^\wedge$ with its hull-kernel topology. The crucial technical tool we use to reach this conclusion is Satake's description of all spherical functions for the pair $(G, U)$, see \cite{satake} and \cite{Mac}. The spherical functions are known to correspond to characters on the abelian
Banach algebra $L^1(G, U)$ associated to the Gelfand pair $(G, U)$. Since our interest is in $*$-representations of $\mathcal{H}(G, U)$, we need to consider a particular class of spherical functions, which we call $*$-spherical functions. We show that the distinction in  isolation of the trivial representation in the natural topologies on $C^*(L^1(G, U))=C^*(p_0L^1(G)p_0)$ and $p_0C^*(G)p_0$, respectively, is due to the fact that not all bounded $^*$-spherical functions on $\mathcal{H}(G, U)$ are positive definite. Indeed, we show that for an arbitrary Gelfand pair $(G, U)$, the canonical surjection from $C^*(L^1(G, U))$ onto $p_0C^*(G)p_0$ is injective precisely when all bounded $*$-spherical functions on $\mathcal{H}(G, U)$ are positive definite, cf. Theorem ~\ref{equivalence spherical functions and canonical Hecke maps I}.

Another consequence of our investigations is that we  can show that for a Gelfand pair $(G, U)$, a necessary and sufficient condition for the existence of the universal $C^*$-completion of $\mathcal{H}(G, U)$ is that all $*$-spherical functions are bounded, see Theorem~\ref{equivalence spherical functions and canonical Hecke maps II}. This behaviour can  perhaps be interpreted as a feature common to the group algebra in the case of abelian $G$, since by taking $U$ to be the trivial subgroup, the spherical functions are characters on $G$, and hence automatically positive definite and bounded.

We thank M. Landstad for indicating that property (T) would be the right notion to consider in order to investigate injectivity of the canonical map when $n\geq 3$. We thank V. Shirbisheh who, after a first version of this paper was circulated, pointed us to references \cite{AD} and \cite{shi}. We thank the anonymous referee for many useful comments which substantially improved the paper and for encouraging us to adopt a more general point of view that encompasses our motivating example: indeed, we are grateful for sharing with us the proof of Theorem~\ref{SLn trivial representation not isolated} in its present generality.

\section{Generalities about spherical functions}

We begin by introducing some notation and recalling terminology. Suppose that $(G, \Gamma)$ is a  Hecke pair, by which we mean either a discrete pair or a pair formed of a topological group with an open subgroup such that every double coset contains finitely many left cosets.  With  $L(g)$ denoting the number of left cosets in $\Gamma g\Gamma$ for $g\in G$, the  map  $\Delta(g)=L(g)/L(g^{-1})$ is a group homomorphism of $G$ into $\mathbb{Q}^+$.  A function $f:G\to \mathbb{C}$ is $\Gamma$-biinvariant if $f(xg)=f(gx)=f(g)$ for all $g\in G$ and $x\in \Gamma$. The \emph{Hecke algebra} $\Hh(G, \Gamma)$ associated to $(G, \Gamma)$ is the $^*$-algebra of $\Gamma$-biinvariant complex valued functions on $G$ which are finitely supported when viewed on $\Gamma\backslash G/\Gamma$; the convolution and involution in $\h(G, \Gamma)$ are defined by
\begin{align}
(f_1\ast f_2)(g)&=\sum_{h\Gamma\in G/\Gamma} f_1(h)f_2(h^{-1}g),\text{ and }\label{eq:conv-Heckealg}\\
f(g)^*&=\Delta(g^{-1})\overline{f(g^{-1}})\label{eq:invol-Heckealg},
\end{align}
for $f_1,f_2, f\in \h(G, \Gamma)$, where the sum in the convolution formula is taken over representatives for the left coset space.

One can define an $L^1$-norm on $\h(G, \Gamma)$ for any Hecke pair $(G, \Gamma)$ by
\begin{equation}\label{eq:L1-norm-Heckealg}
\Vert f\Vert_1=\sum_{\Gamma g\Gamma\in \Gamma\backslash G/\Gamma} \vert f(\Gamma g\Gamma)\vert L(g),
\end{equation}
and take the corresponding completion $L^1(G, \Gamma)$, which is a Banach $^*$-algebra.

An example of a Hecke pair arises when $G$ is a locally compact topological group  and $U$ is a compact open subgroup. This type of example is generic, cf. \cite{tzanev}, see Section~\ref{section:canonical-hom} for more details. Given such $(G, U)$, the Hecke algebra can also be defined as the space of complex valued, compactly supported  functions on $U\backslash G/U$ which are $U$-biinvariant; the convolution and the involution are defined as in \eqref{eq:conv-Heckealg} and \eqref{eq:invol-Heckealg}, where $\Delta$ is  the modular function of $G$.  Let $\mu$ be a left Haar measure on $G$ normalised by $\mu(U)=1$, and consider the Banach $^*$-algebra $L^1(G)$ endowed with its convolution
$$
(f_1\ast f_2)(g)=\int_G f_1(h)f_2(h^{-1}g)\,d\mu(h)
$$
for $f_1,f_2\in L^1(G)$ and involution defined as above. The subspace $L^1(U\backslash G/U)$  of $U$-biinvariant functions in $L^1(G)$ becomes a closed subalgebra of $L^1(G)$. Let $p_0$ be the self-adjoint projection in $L^1(G)$ (and in $C_c(G)$) determined by $\chi_U$. Then $\h(G, U)=p_0C_c(G)p_0$ and  there are canonical isomorphisms $ L^1(U\backslash G/U)\cong L^1(G, U)\cong p_0 L^1(G)p_0$.

Recall  that when $G$ is a \emph{unimodular} locally compact group and $U$ a compact open subgroup, then $(G, U)$ is a \emph{Gelfand pair} if $L^1(U\backslash G/U)$ is commutative. It follows from the results of \cite[\S 2.2]{shi} that whenever $G$ is a locally compact group and $U$ a compact open subgroup such that $\h(G, U)$ is commutative, then $G$ must be unimodular.

The notion of a spherical function for a Gelfand pair is well-known, cf. for example \cite{godement} or the expository article \cite{dieudonne}. Here we follow \cite[\S 5.2]{satake}.

\begin{df}\label{def:spherical-satake} Let $G$ be a unimodular locally compact group and $U$ a compact subgroup.
A continuous function $\omega: G \to \mathbb{C}$ is  a spherical function of $G$ relative to $U$
  if the following conditions are satisfied:
\begin{itemize}
\item[i)] $\omega(e) = 1$.
\item[ii)] $\omega(u_1gu_2) = \omega(g)$ for all $g \in G$ and $u_1, u_2 \in U$.
\item[iii)] $\int_U \omega(g_1u g_2)\,d\mu(u) = \omega(g_1)\omega(g_2)$, for all $g_1, g_2 \in G$.
\end{itemize}
\end{df}

  If $U$ is a compact open subgroup of $G$, then we call a function $\omega$ as in   Definition~\ref{def:spherical-satake}  a  \emph{spherical function} for $\Hh(G, U)$. Assuming $i)$ and $ii)$, condition $iii)$ is equivalent to the following:
\begin{itemize}
\item[iv)] for each $f \in \h(G, U)$ there exists $\lambda_f \in \mathbb{C}$ such that $f*\omega = \lambda_f \omega$.
\end{itemize}

Given a spherical function $\omega$ for $\Hh(G, U)$, one associates a homomorphism $\tau_{\omega}: \h(G, U) \to \mathbb{C}$ defined by
\begin{align}
\tau_{\omega} (f) := \int_G f(g) \omega(g^{-1}) \,d\mu(g)\,.
\end{align}
If $U$ is moreover open, then every non-trivial homomorphism of $\h(G, U)$ into $\mathbb{C}$ arises in this way (a result attributed in \cite{satake} to T. Tamagawa). Thus $\omega \longleftrightarrow \tau_\omega$ establishes a bijective correspondence between spherical functions and homomorphisms of $\h(G, U)$ onto $\mathbb{C}$. This bijective correspondence holds for non-trivial homomorphisms of  $L^1(G, U)$ into  $\mathbb{C}$ when one considers bounded spherical functions. For a proof of the next result, see for example \cite[Theorem 7, Chp. IV]{lang}.

\begin{prop}
\label{bounded spherical functions give L1 characters}
 Let $(G, U)$ be  a Gelfand pair. Then $\omega \longleftrightarrow \tau_\omega$ establishes a bijective correspondence between bounded spherical functions with $|\omega(\cdot)| \leq 1$ and homomorphisms of $L^1(G, U)$ onto $\mathbb{C}$.
\end{prop}

In the next result we identify which spherical functions give rise to $^*$-homomorphisms of $\h(G, U)$. We omit the routine proof.

\begin{prop}
\label{spherical functions and star hom}
 Let $(G, U)$ be a Gelfand pair and $\omega$  a spherical function. The homomorphism $\tau_{\omega}$ is a $^*$-homomorphism if and only if $\omega(g^{-1}) = \overline{\omega(g)}$ for all $g \in G$.
\end{prop}

\begin{df}
 Let $G$ be a unimodular locally compact group and $U$ a compact open subgroup. A spherical function $\omega$ for $\Hh(G, U)$ that satisfies $\omega(g^{-1}) = \overline{\omega(g)}$ for all $g \in G$ will be called a \emph{$^*$-spherical function}.
\end{df}

According to Proposition \ref{spherical functions and star hom}, the $^*$-spherical functions for a Gelfand pair are precisely those for which the associated map $\tau_{\omega}$ is a $^*$-homomorphism.

\begin{prop}
\label{unif convergence on compact sets and weak star top}
 Suppose that $\{\omega_n\}_{n \in \mathbb{N}}$ is a sequence of spherical functions and $\omega$ a spherical function for $\h(G, U)$ such that $|\omega(\cdot)| \leq 1$ and $|\omega_n(\cdot)| \leq 1$ for all $n\geq 1$. If $\omega_{n} \to \omega$ uniformly on compact subsets of $G$, then $\tau_{\omega_n} \to \tau_{\omega}$ in the weak$^*$-topology, as linear functionals of $L^1(G,U)$.
\end{prop}

{\bf \emph{Proof:}} Note that Proposition~\ref{bounded spherical functions give L1 characters} guarantees that $\tau_\omega$ and $\tau_{\omega_n}$, for each $n\geq 1$,  extend to homomorphisms of $L^1(G, U)$ to $\mathbb{C}$. Let $f \in L^1(G,U)$ and $\epsilon >0$. Since $f$ can be viewed as a function in $L^1(G)$, we can choose a sufficiently large compact set $K \subseteq G$ such that $\int_{G \backslash K} f \;d\mu < \epsilon / 4$.  Choose $n_0 \in \mathbb{N}$ such that for $n \geq n_0$ we have $
\|(\omega_n - \omega)|_{K^{-1}} \|_{\infty} < {\epsilon}/{2 \|f\|_1}\,.$
Then for $n \geq n_0$ we have
\begin{eqnarray*}
| \tau_{\omega_n}(f) - \tau_{\omega}(f)| & \leq & \int_{G \backslash K} |f(g)||\omega_n(g^{-1}) - \omega(g^{-1})| \;d\mu(g)\\
& + & \int_{K} |f(g)|\;d\mu(g) \cdot \|(\omega_n - \omega)|_{K^{-1}} \|_{\infty}\\
& \leq & \frac{2 \epsilon}{4} + \frac{\epsilon}{2}=  \epsilon\,.
\end{eqnarray*}
{}

Given a Gelfand pair $(G, U)$, a spherical function $\omega:G\to \CC$ is \emph{positive definite} if for any $n\geq 1$, any finite subset $\{s_1,\dots ,s_m\}$ of $G$ and any  collection $\{z_1,\dots, z_m\}$ of complex numbers, one has
$$
\sum_{j,k=1}^m \omega (s_j^{-1}s_k)\overline{z_j}z_k\geq 0,
$$
see for example \cite{dieudonne}. A slightly different formulation of this notion appears in \cite[Appendix II]{satake}.
Positive definite spherical functions for a Gelfand pair $(G,U)$ appear naturally within the context of unitary representations. In the general theory of Hecke algebras a unitary representation of $G$ on a Hilbert space which is generated by its $U$-fixed vectors always gives rise to a nondegenerate $^*$-representation of the Hecke algebra $\h(G,U)$. The converse, whether every nondegenerate $^*$-representation of $\h(G,U)$ arises in this way from a unitary representation, is not always true. Hall proposed a positivity condition for a certain inner product as a sufficient condition for a category equivalence of the two classes of representations, \cite[Theorem 3.25]{hall}. The category equivalence holds true when certain $C^*$-completions of the Hecke algebra are the same, see \cite[Corollaries 6.11(ii) and 6.19]{kal-land-qui}. In the case of Gelfand pairs, bounded positive definite spherical functions $\omega$ are precisely those for which the associated $^*$-representation $\tau_{\omega}$ arises from a unitary representation of $G$. Before we indicate the relation between positive definite spherical functions and $C^*$-completions of Hecke algebras, we first recall  that positive definite spherical functions are always bounded and $^*$-spherical. The result is well-known, cf. for example \cite[Appendix II]{satake} or \cite[\S5, Chp. IV]{lang}.

\begin{prop}
Let $(G,U)$ be a Gelfand pair. Every positive definite spherical function is necessarily bounded and $^*$-spherical.
\end{prop}

\section{The canonical surjection for a Hecke pair and $*$-spherical functions}\label{section:canonical-hom}

To motivate our first result, we recall that for a (discrete) Hecke pair $(G, \Gamma)$, the convolution and involution
on the associated Hecke algebra $\Hh(G, \Gamma)$ are as given in equations \eqref{eq:conv-Heckealg} and \eqref{eq:invol-Heckealg}. It was shown in \cite{tzanev}, see also \cite{kal-land-qui}, that there is a unique pair $(\overline{G}, \overline{\Gamma})$ where $\overline{G}$ is a totally disconnected locally compact group, $\overline{\Gamma}$ is a compact open subgroup, there is a canonical  embedding $\iota:G\to \overline{G}$ such that $\iota(G)$ is dense in $\overline{G}$, $\iota(\Gamma)$ is dense in $\overline{\Gamma}$, and $\iota^{-1}(\overline{\Gamma})=\Gamma$ (more precisely, the uniqueness is achieved by passing to the reduction $(G_r, \Gamma_r)$ of $(G, \Gamma)$).  We refer to $(\overline{G}, \overline{\Gamma})$ as the Schlichting completion of $(G, \Gamma)$.

Let $p_0$ denote the self-adjoint projection $\chi_{\overline{\Gamma}}$  in $C_c(\overline{G})$. There are canonical isomorphisms $
\Hh(G, \Gamma)\cong \Hh(\overline{G}, \overline{\Gamma})\cong p_0 C_c(\overline{G})p_0$. Upon completion with respect to the norm from \eqref{eq:L1-norm-Heckealg}, one obtains an isomorphism $L^1(G, \Gamma)\cong p_0L^1({\overline{G}})p_0$.

 The canonical surjection for the Hecke pair $(G, \Gamma)$ alluded to in the title of the section is the $*$-homomorphism
 \begin{equation}\label{eq:canonical-surj}
 \Pi:C^*(L^1(G, \Gamma))\longrightarrow p_0C^*(\overline{G})p_0,
 \end{equation}
 see \cite{tzanev} and \cite{kal-land-qui}.
 This map is known to be an isomorphism in many cases, such as when $\overline{G}$ is Hermitian, \cite{kal-land-qui}, or when $L^1(\overline{G})$ is quasi-symmetric in the terminology of \cite{palma1}. The last property is known to hold for a class of groups containing Hermitian groups and groups with subexponential growth. For a Gelfand pair, the next result gives a necessary and sufficient condition for the canonical surjection to be an isomorphism in terms of properties of $*$-spherical functions.

\begin{thm}
\label{equivalence spherical functions and canonical Hecke maps I}
 Let $G$ be a unimodular locally compact group and $U$ a compact open subgroup such that $(G, U)$ is a Gelfand pair. The following are equivalent:

\textnormal{(a)} All bounded $^*$-spherical functions for $\h(G, U)$ are positive definite.

\textnormal{(b)}  The canonical surjection $\Pi: C^*(L^1(G, U))\longrightarrow p_0C^*(G)p_0$
is an isomorphism.
\end{thm}

{\bf \emph{Proof:}}
We have that $p_0$ is the self-adjoint projection in $L^1(G)$ equal to $\chi_U$.
We aim to use \cite[Corollary 6.11(ii)]{kal-land-qui}, according to which
the Banach $^*$-algebra $D:= p_0L^1(G)p_0$ is such that $C^*(D)=p_0C^*(G)p_0$ if and only if every $^*$-representation $\pi$ of $D$ on a Hilbert space satisfies
\begin{equation}\label{eq:right-ip-positive}
\pi(\langle f, f\rangle_D)\geq 0 \text{ for all }f\in L^1(G)p_0,
\end{equation}
where the $D$-valued inner product $\langle \cdot, \cdot \rangle_D$ is defined by $\langle f_1, f_2\rangle_D=f_1^*f_2$ for $f_1, f_2\in L^1(G)p_0$.

Let us assume (a). We will establish \eqref{eq:right-ip-positive}, which therefore shows that $C^*(D)=p_0C^*(G)p_0$.

{\bf{Claim 1}}: The inequality in \eqref{eq:right-ip-positive} is valid when $\pi$ runs over the set of characters of $D$. To see this,  let $f=\sum_{j=1}^m z_js_j p_0$ in $C_c(G)p_0\subset L^1(G)p_0$; then $f^*=\sum_{j=1}^m \overline{z_j}p_0s_{j}^{-1}$. Let $\omega$ be a bounded $^*$-spherical function, which by assumption is positive definite. Then, using that $p_0gp_0=\frac 1{L(g)} \chi_{UgU}$ for every $g\in G$, it follows that
\begin{align}
\tau_\omega(\langle f, f\rangle_D)
&= \tau_\omega (p_0\sum_{j,k=1}^m \overline{z_j}z_ks_j^{-1}s_k p_0)\notag \\
&= \sum_{j,k=1}^m \overline{z_j}z_k\tau_\omega(p_0 s_j^{-1}s_k p_0)\notag \\
&= \sum_{j,k=1}^m \frac {\overline{z_j}z_k}{L(s_j^{-1}s_k)}\int_G \chi_{Us_j^{-1}s_kU}(g)\omega(g^{-1})d\mu(g)\notag \\
&= \sum_{j,k=1}^m \frac {\overline{z_j}z_k}{L(s_j^{-1}s_k)}\int_{Us_j^{-1}s_kU} \omega(s_j^{-1}s_k)d\mu(g)\notag \\
&=\sum_{j,k=1}^m  \omega(s_j^{-1}s_k)\overline{z_j}z_k.\label{tau-omega-on-inner-product}
\end{align}
Thus $\tau_\omega(\langle f, f\rangle_D)\geq 0$ for all bounded $^*$-spherical functions $\omega$ on $\h(G, U)$. By continuity of $\tau_\omega$, it follows that $\tau_\omega(\langle f, f\rangle_D)\geq 0$ for all $f\in L^1(G)p_0$ and all bounded $^*$-spherical functions on $\h(G, U)$. Since $D\cong L^1(G, U)$ and every character of $L ^1(G, U)$ is of the form $\tau_\omega$ for some bounded $^*$-spherical function $\omega$, it follows that
$$
\pi(\langle f, f\rangle_D)\geq 0
$$
for any character $\pi$ of $L^1(G, U)$. This proves claim 1.

{\bf{Claim 2}}: If \eqref{eq:right-ip-positive} is valid when $\pi$ runs over the set of characters of $D$, then it is valid for arbitrary $\pi$.
 To show this, let $\pi$ be a $^*$-representation of the commutative, unital Banach $^*$-algebra $D$ on a Hilbert space $H_\pi$ and fix an arbitrary state $\psi$ of $B(H_\pi)$. Then $\psi\circ \pi$ is a state on $D$ is the sense of $*$-algebras, see for example \cite[Definition 9.4.21]{palmer}. Therefore, by \cite[Theorem 9.6.6]{palmer}, we have $\psi\circ \pi\in \overline{\operatorname{conv}}\{D_P\}$, where $D_P$ denotes the set of pure states on $D$. Since $D$ is a commutative $*$-algebra,  \cite[Theorem 9.6.10]{palmer} says that the set $D_P$ of pure states is the same as the set of $*$-homomorphisms of $D$ onto $\mathbb{C}$. The latter one is the character space of $D$. Now, $\psi\circ \pi$ admits a barycentric decomposition
 $$
 (\psi\circ \pi)(d)=\int \tau'(d)d\mu(\tau'),
 $$
 for a unique measure $\mu$ supported on the set of pure states (which are the extremal points in the set of states on $D$), see \cite[\S 4.1]{bra-rob}. Since $\tau'(\langle f, f\rangle_D)\geq 0$ for every $f\in L^1(G)p_0$ by claim 1, it follows that $(\psi\circ \pi)(\langle f, f\rangle_D)\geq 0$ for all $f\in L^1(G)p_0$. Since $\psi$ was chosen arbitrarily, it follows that $\pi(\langle f, f\rangle_D)\geq 0$, proving claim 2 and finishing the proof of one implication.

Conversely, assume (b). Let $\omega$ be a bounded $^*$-spherical function on $\h(G, U)$. By Propositions~\ref{bounded spherical functions give L1 characters} and \ref{spherical functions and star hom}, $\tau_\omega$ extends to a $^*$-homomorphism  $L^1(G, U)\to \CC$. Fix a finite subset $\{s_1,\dots ,s_m\}$ of $G$ and a finite subset $\{z_1,\dots, z_m\}$ of complex numbers. Let $f=\sum_{j=1}^m z_js_jp_0\in L^1(G)p_0$. Since $\tau_\omega$ satisfies
\eqref{eq:right-ip-positive},  a computation similar to the one leading to equation \eqref{tau-omega-on-inner-product} shows that
$\sum_{j,k=1}^m \omega (s_j^{-1}s_k)\overline{z_j}z_k=\tau_\omega(\langle f, f\rangle_D)\geq 0$, as needed in (a).
{}\\

The Hecke algebra $\Hh(G, \Gamma)$ associated to a Hecke pair $(G, \Gamma)$ need not admit a universal $C^*$-completion, \cite{hall}. When it does, the universal $C^*$-completion $C^*(G, \Gamma)$ of $\Hh(G, \Gamma)$ admits a natural surjection onto $C^*(L^1(G,\Gamma))$. As pointed out in \cite[Questions 6.16(ii)]{kal-land-qui}, this map is an isomorphism for all known classes such that $C^*(G, \Gamma)$ exists,  although a general explanation of why this must be the case is missing. The next result gives a necessary and sufficient condition for existence of a universal $C^*$-completion of a Gelfand pair. As  in all other similar cases, the natural surjection then becomes injective.

\begin{thm}
\label{equivalence spherical functions and canonical Hecke maps II}
 Let $G$ be a unimodular locally compact group and $U$ a compact open subgroup such that $(G, U)$ is a Gelfand pair. The following are equivalent:

\textnormal{(a)} All $^*$-spherical functions for $\h(G, U)$ are bounded.

\textnormal{(b)} The universal $C^*$-completion $C^*(G, U)$ exists and the canonical surjection $C^*(G, U)\longrightarrow C^*(L^1(G, U))$ is an isomorphism.
\end{thm}

{\bf \emph{Proof:}} Let us assume (a) is true. Since $\h(G, U)$ is abelian, in order to prove that $C^*(G, U)$ exists it suffices to show that $\sup_{\phi} |\phi(f)| < \infty$ for every $f \in \h(G,U)$, where the supremum runs over the set of $^*$-homomorphisms $\phi:\h(G, U) \to \mathbb{C}$. We know that each $^*$-homomorphism $\h(G, U) \to \mathbb{C}$ is of the form $\tau_{\omega}$ for some $^*$-spherical function $\omega$. By (a) every $^*$-spherical function is bounded, and therefore $\tau_{\omega}$ extends to a character of $L^1(G,U)$. Hence
\begin{align*}
|\tau_{\omega}(f)| \leq \|f\|_{L^1(G,U)}\,,
\end{align*}
from which it follows immediately that $\sup_{\phi} |\phi(f)|\leq \|f\|_{L^1(G,U)}  < \infty$. We conclude that $C^*(G, U)$ exists. Moreover, since every $^*$-homomorphism $\h(G,U) \to \mathbb{C}$ extends to a character of $L^1(G,U)$ it follows that the canonical surjection $C^*(G, U)\longrightarrow C^*(L^1(G, U))$ is an isomorphism.

Let us now assume (b). Let $\omega$ be a $^*$-spherical function. As we know, $\tau_{\omega}$ is then a $^*$-homomorphism $\h(G,U) \to \mathbb{C}$. By (b) we must have
\begin{align*}
 |\tau_{\omega}(f)| \leq \|f\|_{C^*(G,U)} = \| f\|_{C^*(L^1(G,U))} \leq \| f \|_{L^1(G,U)}\,.
\end{align*}
Hence, $\tau_{\omega}$ extends to a character of $L^1(G,U)$. Thus, $\omega$ must be bounded. {}\\

An immediate consequence of Theorems~\ref{equivalence spherical functions and canonical Hecke maps I} and \ref{equivalence spherical functions and canonical Hecke maps II} is the following:

\begin{cor}
\label{equivalence spherical functions and canonical Hecke maps III}
 Let $G$ be a unimodular locally compact group and $U$ a compact open subgroup such that $(G, U)$ is a Gelfand pair. The following are equivalent:

\textnormal{(a)} All $^*$-spherical functions for $\h(G, U)$ are bounded and positive definite.

\textnormal{(b)} the universal $C^*$-completion $C^*(G, U)$ exists and the natural maps are isomorphisms:
$$
C^*(G, U)\overset{\cong}{\longrightarrow} C^*(L^1(G, U))\overset{\cong}{\longrightarrow} p_0C^*(G)p_0.
$$
\end{cor}

\section{Spherical functions for simple algebraic groups over $\mathfrak{p}$-adic fields}

The spherical functions for the pair  $(SL_n(\mathbb{Q}_p), SL_n(\mathbb{Z}_p))$ and, more generally, reductive algebraic groups over $\mathfrak{p}$-adic fields, were completely characterised by Satake, see \cite{satake-short} and \cite{satake}. The results in \cite{satake} are valid in greater generality, see for example \cite{Mac}.

We introduce some notation first: $\mathfrak{K}$ denotes a finite extension of $\mathbb{Q}_p$, $\mathfrak{o}$ is the ring of integers in $\mathfrak{K}$, $\pi$ is a prime element generating the unique prime ideal $\mathfrak{p}$ of $\mathfrak{o}$, and $q$ is the number of elements in the residue class field $\mathfrak{o}/\mathfrak{p}$. An algebraic group over $\mathfrak{K}$ is the group of $\mathfrak{K}$-rational points of an algebraic group defined over $\mathfrak{K}$. Our interest is in simple algebraic groups over $\mathfrak{K}$ of rank $\nu$, with $\nu\geq 2$.

The main objects in \cite{satake} are $G$, $U$, $H$, $N$, where $G$ is an algebraic group over $K$,  $U$ is a ``good'' open compact subgroup of $G$, $H$ is a closed subgroup consisting of semi-simple elements and $N$ is a unipotent subgroup of $G$ normalised by $H$ such that if $H^u$ is the unique maximal compact subgroup of $H$, the following condition (I) is satisfied: $G=UHN$ and $U\supset H^u$. In addition, the Weyl group of $G$ relative to $H$, $W_H$, must satisfy a technical condition (II) which depends on the existence of a fundamental domain of $W_H$ in $H/H^u$, see \cite[\S 3.3]{satake}. Then $(G, U)$ is a Gelfand pair, and under the assumptions (I) and (II), \cite[Theorem 2, \S 5.4]{satake} gives a parametrisation of all the spherical functions for $\h(G, U)$. Given a semi-simple algebraic group $G$ over $K$, subgroups $U$, $H$, $N$ satisfying the assumptions (I) and (II) can be
found in many relevant examples, see \cite[Chapter III]{satake} and \cite{Mac}.

We assume for the remaining of this section that $G$, $U$, $H$, $N$ satisfy the assumptions (I) and (II). We need to recall the parametrisation of spherical functions from \cite[Theorem 2, \S 5.4]{satake}, see also \cite[Theorem 1]{Mac}.

Following \cite{satake}, we let $X(H)$ be the set of $\mathfrak{K}$-morphisms of $H$ into $\mathfrak{K}^*$ (this is the character group of H). The subgroup $H^u=\{h\in H\mid \chi(h)\in \mathfrak{o}^*\text{ for all }\chi\in X(H)\}$ is the unique maximal compact subgroup of $H$. Letting $M\subset \operatorname{Hom}(X(H), \mathbb{Z})$ be the subgroup formed of elements $l_h:X(H)\to \mathbb{Z}$ given by $l_h(\chi)=\operatorname{ord}_{\mathfrak{p}}(\chi(h))$ for $\chi \in X(H)$ and $h\in H$, the map $h\mapsto l_h$ induces an isomorphism
\begin{equation}\label{eq:iso-HoverHu-M}
H/H^u\cong M.
\end{equation}
We let $\pi^{\mathbf{m}}$ denote the preimage of $\mathbf{m}\in M$ under this isomorphism.
The Weyl group $W=W_H$ acts on $M$ by (inner) automorphisms. There is a natural pairing
\begin{equation}\label{eq:pairing}
 M \times (X(H) \otimes \mathbb{C})\mapsto \mathbb{C},\, (\mathbf{m}, \mathbf{s})\mapsto \mathbf{m}.\mathbf{s},
\end{equation}
and under this pairing $W$ acts on $X(H) \otimes \mathbb{C}$.

\begin{rem} In the case $G=SL_n(\mathfrak{K})$, then $U=SL_n(\mathfrak{o})$, $H$ is the subgroup of diagonal matrices, $H^u$ is the subgroup of diagonal matrices $h = \mathrm{diag}(h_1, \dots, h_n)$ in $H$ such that $h_k \in \mathfrak{o}^*$ for $k=1,\dots, n$, $M\cong \{m\in \mathbb{Z}^n\mid \sum_{k = 1}^n m_k = 0\}$, and the Weyl group $W$ is the symmetric group $S_n$ on $n$ letters acting on $M$ by permutations. If $\mathbf{m}= (m_1, \dots, m_n)\in M$, then $\pi^{\mathbf{m}} := \mathrm{diag}(\pi^{m_1}, \dots, \pi^{m_n})\in H$. There is a natural identification  $X(H) \otimes \mathbb{C} \cong \mathbb{C}^n / \{(s, \dots, s): s \in \mathbb{C} \}$ under which the pairing \eqref{eq:pairing} takes the form $\mathbf{m}.\mathbf{s} = \sum_{k = 1}^n m_ks_k$,
for $\mathbf{m} = (m_1, \dots, m_n) \in M$ and $\mathbf{s} = (s_1, \dots, s_n) \in \mathbb{C}^n$ representing a class in $X(H) \otimes \mathbb{C}$.
\end{rem}

Let $
 \widehat{M}: = \{\mathbf{s} \in X(H) \otimes \mathbb{C} : \mathbf{m}.\mathbf{s}\in \mathbb{Z} \;\;\text{for all} \;\; \mathbf{m} \in M \}$. By \cite[\S 5.4]{satake}, if $\alpha:H \to \mathbb{C}^*$ is a quasi-character (i.e. a continuous homomorphism of $H$ into $\mathbb{C}^*$ with respect to the $\mathfrak{p}$-adic topology) satisfying $\alpha(H^u) = 1$, then $\alpha$ is uniquely determined by the values of $\alpha(\pi^{\mathbf{m}})$ with $\mathbf{m} \in M$. Moreover, there exists $\mathbf{s} \in X(H) \otimes \mathbb{C}$, uniquely determined modulo $ \frac{2\pi i}{\log q} \widehat{M}$,  such that
\begin{align}
\label{eq for alpha with respect to s}
\alpha(\pi^{\mathbf{m}}) = q^{- \mathbf{m}.\mathbf{s}}\,.
\end{align}
%Here, if $\langle \cdot, \cdot \rangle$ denotes the natural $\mathbb{Z}$-valued pairing between $X(H)$ and $\operatorname{Hom}(X(H),\mathbb{Z})$, we %identify $\mathbf{m}.\mathbf{s}$ as $-\bigl((2\pi i)/{\log q}\bigr) \langle \mathbf{m}, \mathbf{s}\rangle$.
Whenever $\alpha(H^u) = 1$ and (\ref{eq for alpha with respect to s}) holds, we use Satake's notation and write $\alpha \leftrightarrow \mathbf{s}$.

Let $\delta$ be the positive quasi-character introduced in \cite[\S 4.1]{satake}; thus, if $d$ is the Haar measure on $N$ normalised so that $\int_{N\cap U}dn=1$, then
$$
\delta(h)=\frac {d(hnh^{-1})}{d(n)} \text{ for }h\in H, n\in N.
$$
Let now $\delta^{\frac{1}{2}}\alpha\leftrightarrow \mathbf{s}$. For $g\in G$ of the form $g=uhn$ with $u\in U$, $h\in H$ and $n\in N$, define a complex-valued function $\psi_{\alpha}$ on $G$ by $\psi_{\alpha} (uhn) := \alpha(h)$. Let $du$ be normalised Haar measure on $U$. Then, according to \cite[\S 5.3 and \S 5.4]{satake}, the formula
\begin{align}\label{eq:psi-from-alpha}
\omega_{\mathbf{s}}(g) := \int_{U} \psi_{\alpha} (g^{-1} u)\, du
\end{align}
for $g\in G$ determines a spherical function $\omega_{\mathbf{s}}$ for $\h(G, U)$.

Before stating the next result we make the following observation: if $\delta^{\frac{1}{2}}\alpha\leftrightarrow \mathbf{s}$, let $\overline{\mathbf{s}}$ be given by $\alpha_{\overline{\mathbf{s}}}(h)=\overline{\alpha_{\mathbf{s}}(h)}$ for all $h\in H$. Then for each $g\in G$ we have
$$
\overline{\omega_{\mathbf{s}}(g^{-1})}=\int_U \overline{\psi_{\alpha} (g u)}\, du
=\int_U {\psi_{\overline{\alpha}} (g u)}\, du=\omega_{\overline{\mathbf{s}}}(g^{-1}).
$$

\begin{lemma}\label{lemma:cond-star-spherical}
If  $\mathbf{s}=-w\overline{\mathbf{s}}$ for some $w$ in $W$, then $\omega_{\mathbf{s}}$ is a $*$-spherical function.
\end{lemma}

{\bf \emph{Proof:}}
To see this, fix $s$ and $w$  with $\mathbf{s}=-w\overline{\mathbf{s}}$, and recall from \cite[Proposition 5.2]{satake} that $\omega_{w\mathbf{s}}=\omega_{\mathbf{s}}$ and $\omega_{-\mathbf{s}}(g) = \omega_{\mathbf{s}}(g^{-1})$ for all $g\in G$.  Then for each $g\in G$ we have $\omega_{\mathbf{s}}(g)=\omega_{w^{-1}\mathbf{s}}(g)=\omega_{-\overline{\mathbf{s}}}(g)=\omega_{\overline{\mathbf{s}}}(g^{-1})$. By the observation preceding the lemma, it then follows that $\overline{\omega_{\mathbf{s}}(g)}={\omega_{\mathbf{s}}(g^{-1})}$, as claimed.{}\\

Now,  \cite[Theorem 2, \S 5.4]{satake} or \cite[Theorem 1]{Mac} establish that every spherical function for $\h(G, U)$ arises as in \eqref{eq:psi-from-alpha} from some $\alpha$ (or, equivalently, from $\mathbf{s}$ such that $\delta^{\frac{1}{2}}\alpha\leftrightarrow \mathbf{s}$). Moreover, $\omega_{\mathbf{s}}=\omega_{\mathbf{s}'}$ if and only if ${\mathbf{s}'}={w\mathbf{s}}$ for some $w\in W$, where the action of $W$ on $M$ is carried over to an action of $W$ by $\CC$-linear transformations of $X(H)\otimes \CC$ which leave $\widehat{M}$ invariant.

 It is immediate to see that the trivial spherical function $\omega_{\mathbf{t}}\equiv 1$ is attained at (the orbit of) $\mathbf{t}\in \bigl(X(H) \otimes \mathbb{C}\bigr)/\bigl(\frac{2\pi i}{\log q} \widehat{M}\bigr)$ satisfying that $\delta(\pi^{\mathbf{m}})^{\frac 12}=q^{\mathbf{m}. \mathbf{t}}$ for all $\mathbf{m}\in M$.

Our initial approach to the following theorem covered the case of  $(SL_n(\mathfrak{K}), SL_n(\mathfrak{o}))$ for $n\geq 3$. We are grateful to the anonymous referee who generously shared the proof in the case of an arbitrary pair $(G, U)$ satisfying Satake's conditions (I) and (II).

\begin{thm}
\label{spherical functions approximating trivial spherical function}
 Let $G, U, H, N$ be groups satisfying conditions (I) and (II) of Satake. Let $(r_n)_{n\geq 1}$ be a non-decreasing sequence of positive real numbers that converges to $1$.  Then  $\{\omega_{r_n\mathbf{t}} \}_{n \in \mathbb{N}}$ is a sequence of bounded $^*$-spherical functions converging to $1$ uniformly on compact subsets of $G$. Moreover, the functions $\omega_{r_n\mathbf{t}}$ are pairwise distinct for large enough $n$.\\
\end{thm}

{\bf \emph{Proof:}} We show first that each $\omega_{r_n\mathbf{t}}$ is $*$-spherical. Note that $\overline{\mathbf{t}}=\mathbf{t}$. The assumption
$\delta^{1/2}\leftrightarrow \mathbf{t}$ implies that $\mathbf{t}$ corresponds to half the sum of the positive roots in the reduced root system of $(G, H)$. Hence by \cite[Chap. VI, \S1, Corollary 3]{bour}, there is $w\in W$ such that $w\mathbf{t}=-\mathbf{t}$. Since $r\mathbf{t}=\overline{r\mathbf{t}}$ for every positive real number $r$, it follows from Lemma~\ref{lemma:cond-star-spherical} that $\omega_{r_n\mathbf{t}}$ is $*$-spherical for each $n\geq 1$.

 To see that the functions $\omega_{r_n\mathbf{t}}$ are pairwise distinct for large  $n$ we note that the affine Weyl group, which is given by
 $W_a=W. \bigl(\frac{2\pi i}{\log q} \widehat{M}\bigr)$, acts properly on $X(H)\otimes \mathbb{C}$, see \cite[Chap. VI, \S2, Proposition 2]{bour}. Hence, for large enough $n\neq m$, the elements $r_n\mathbf{t}$ and $r_m\mathbf{t}$ in  $X(H)\otimes \mathbb{C}$ belong to different $W_a$ orbits. Then by \cite[Chap. II, \S 5, Theorem 2]{satake}, the spherical functions $\omega_{r_n\mathbf{t}}$ and $\omega_{r_m\mathbf{t}}$ are distinct.

 We next turn to proving that $\omega_{r_n\mathbf{t}}$ is bounded by $1$ for each $n\geq 1$. For $z\in \mathbb{C}$, note that $\delta^{z/2}\leftrightarrow z\mathbf{t}$. We claim that $|\omega_{z\mathbf{t}}(\cdot)|\leq 1$ when $|Re(z)|\leq 1$. For a fixed $g\in G$, the function $z\mapsto \omega_{z\mathbf{t}}(g)$ is entire. For a given $u\in U$ we write $g^{-1}u=u'h_un'$ according to the decomposition $G=UHN$. By \eqref{eq:psi-from-alpha}, $ \omega_{z\mathbf{t}}(g)=\int_U \delta(h_u)^{(z-1)/2}du$. Thus, $|\omega_{z\mathbf{t}}(g)|\leq 1$ whenever $Re(z)=1$. If $Re(z)=-1$, using that $\omega_{z\mathbf{t}}=\omega_{-z\mathbf{t}}(g^{-1})$, it again follows that $|\omega_{z\mathbf{t}}(g)|\leq 1$. For a general $z \in \mathbb{C}$ we have that
\begin{align*}
|\omega_{z\mathbf{t}}(g)| \leq \sup_{u \in U} |\delta(h_u)^{(z-1)/2}| = \sup_{u \in U} |\delta(h_u)|^{(Re(z)-1)/2}\,,
\end{align*}
and therefore there is a constant $C$ such that $|\omega_{z\mathbf{t}}(g)| < C$ for any $z$ with $|Re(z)|\leq 1$.
It now follows from Hadamard's three-lines theorem, applied to the function $z \mapsto \omega_{z\mathbf{t}}(g)$ on the strip $|Re(z)| \leq 1$, that $|\omega_{z\mathbf{t}}(g)| \leq 1$ for any $z$ on the strip. Since this is true for any $g \in G$, the claim follows.

It remains to show that $\omega_{r_n\mathbf{t}} \to \omega_{\mathbf{t}}$ uniformly on compact sets. Since spherical functions are constant on the compact open sets $UgU$, with $g \in G$, it is enough to prove pointwise convergence, i.e. $\omega_{r_n\mathbf{t}}(g) \to \omega_{\mathbf{t}} (g)$ for any $g\in G$. Now, if $\delta^{1/2}\alpha_n\leftrightarrow r_n\mathbf{t}$, then $\alpha_n$ converges to $1$ in the parameter space for spherical functions. Passing through the identity \eqref{eq:psi-from-alpha}, the claim follows from a routine argument using the Lebesgue dominated convergence theorem.

\section{Property (T) and the canonical surjection for $(G, U)$}

Let $A$ be a $^*$-algebra.  We will denote by $\mathrm{Prim}(A)$ the \emph{primitive ideal space} of $A$, i.e. the set of kernels of topologically irreducible $^*$-representations of $A$, with the hull-kernel topology. Moreover, we will denote by $\widehat{A}$  the \emph{dual space} of $A$, which is the set of unitary equivalence classes of irreducible $^*$-representations of $A$. There is a canonical map $\widehat{A} \to \mathrm{Prim}(A)$ defined by $[\pi] \mapsto Ker(\pi)$, and the topology of $\widehat{A}$ is the topology pulled back from $\mathrm{Prim}(A)$ through this map.

Property (T) for a topological group is known to have several equivalent formulations. For our purposes two of these characterisations are crucial: one of them because it captures a property (T) for Hecke pairs, as introduced by Tzanev in \cite{tzanev-thesis}, see Appendix A and in particular Theorem~\ref{thm:equiv-T-Heckepair}. The other characterisation has our interest because it \emph{fails} to pass to Hecke pairs in the sense of Theorem~\ref{SLn trivial representation not isolated}. As we shall now show, this phenomenon is intimately related to the question of whether the canonical surjection from \eqref{eq:canonical-surj} is injective.

\begin{prop}\label{prop:isolation-passes-to-corner}
Let $(G, \Gamma)$ be a Hecke pair with Schlichting completion $(\overline{G}, \overline{\Gamma})$. Suppose that $\overline{G}$ has property (T) or, equivalently, by Theorem~\ref{thm:equiv-T-Heckepair} the pair $(G, \Gamma)$ has property (T) for Hecke pairs.
Then the trivial representation of $p_0C^*(\overline{G})p_0$ is an isolated point of $\big(p_0C^*(\overline{G})p_0 \big)^{\wedge}$ with the hull-kernel topology.
\end{prop}

{\bf \emph{Proof:}} Since $\overline{G}$ has property (T), the trivial representation $\mathbf{1}_{\overline{G}}$ of $\overline{G}$ is isolated  in $\widehat{\overline{G}}$ with the Fell topology. Under the standard identifications, the Fell topology on $\widehat{\overline{G}}$ coincides with the hull-kernel topology on $C^*(\overline{G})^{\wedge}$ (see \cite[\S 18]{dixmier} and comment on page 413 of [BdHV]), so that the trivial representation of $C^*(\overline{G})$ is an isolated point of $\big(C^*(\overline{G}) \big)^{\wedge}$. Moreover, by \cite[Proposition A.27 (b)]{raeburn williams}, we  see that $\mathbf{1}_{\overline{G}}$, restricted to the ideal $\overline{C^*(\overline{G})p_0C^*(\overline{G})}$, is an isolated point of $\Big(\overline{C^*(\overline{G})p_0C^*(\overline{G})} \Big)^{\wedge}$. By the Morita equivalence between $\overline{C^*(\overline{G})p_0C^*(\overline{G})}$ and $p_0C^*(\overline{G})p_0$, together with \cite[Corollary 3.33 (a)]{raeburn williams}, we see that the trivial representation of $p_0C^*(\overline{G})p_0$ is an isolated point of $\big(p_0C^*(\overline{G})p_0 \big)^{\wedge}$. {}

The above proposition shows that $p_0C^*(\overline{G})p_0$ still captures the characterisation of property (T) in terms of the isolation of the trivial representation. However, as the next result shows, this characterisation of property (T) fails in general to pass to $C^*(L^1(\overline{G}, \overline{\Gamma}))$.

\begin{thm}
\label{SLn trivial representation not isolated}
Let $(G, U)$ be a Gelfand pair as in Theorem~\ref{spherical functions approximating trivial spherical function}. The trivial representation of $C^*(L^1(G,U))$ is not an isolated point of the space $\big(C^*(L^1(G,U)) \big)^{\wedge}$ endowed with its  hull-kernel topology.
\end{thm}

{\bf \emph{Proof:}} By Theorem \ref{spherical functions approximating trivial spherical function} there is a sequence of mutually different bounded $^*$-spherical functions $\{\omega_{\mathbf{t}_n}\}_{n \in \mathbb{N}}$ that converges uniformly on compact sets to the trivial spherical function  $\omega_{\mathbf{t}}$.

By Proposition \ref{bounded spherical functions give L1 characters}, each $\tau_{\omega_{\mathbf{t}_n}}$ is a $^*$-homomorphism of $L^1(G,U)$ onto $\mathbb{C}$. Moreover, by Proposition \ref{unif convergence on compact sets and weak star top}, the sequence $\{\tau_{\omega_{\mathbf{t}_n}}\}_{n\in \mathbb{N}}$ converges to the trivial representation $\tau_{\omega_{\mathbf{t}}}$ of $L^1(G,U)$ in the weak$^*$-topology.

Let us denote by $L^1(G, U)^{\wedge}$ the space of non-trivial $^*$-homomorphisms of $L^1(G,U)$ onto $\mathbb{C}$ endowed with the weak$^*$-topology, which is the same as  the pure state space of $L^1(G,U)$, since this Banach $^*$-algebra is abelian. Since the $^*$-homomorphisms $\tau_{\omega_{\mathbf{t}_n}}$ are mutually different, the trivial representation of $L^1(G,U)$ is not an isolated point of $L^1(G, U)^{\wedge}$.

By \cite[Theorem 10.1.12 (b)]{palmer} the spaces $\big(C^*(L^1(G,U)) \big)^{\wedge}$ and $L^1(G,U)^{\wedge}$ are homeomorphic. Hence, the trivial representation of $C^*(L^1(G,U))$ is not an isolated point of $\big(C^*(L^1(G,U)) \big)^{\wedge}$ with the weak$^*$-topology. Since $C^*(L^1(G,U))$ is abelian, the weak$^*$-topology and the hull-kernel topologies coincide (by, for example, \cite[Appendix A]{raeburn williams}), and this proves our claim. {}\\

 It is known that a simple algebraic group  of rank $\nu\geq 2$ over a $\mathfrak{p}$-adic field has property (T), see for example \cite[Theorem 1.6.1]{BdHV}. Hence by Proposition~\ref{G has T iff G with compact has T}, $(G, U)$ has property (T) for Hecke pairs. Further, by Proposition~\ref{prop:isolation-passes-to-corner} we have  that the trivial representation of $p_0C^*({G})p_0$  is isolated in $\big(p_0C^*({G})p_0\big)^{\wedge}$. This fact and Theorem~\ref{SLn trivial representation not isolated} lead to the following corollary:

\begin{cor}
 Let $(G,U)$ be a Gelfand pair as in Theorem~\ref{SLn trivial representation not isolated} and assume that $G$ is simple and of rank $\nu\geq 2$. Then the $C^*$-completions $C^*(L^1(G, U))$ and $p_0C^*({G})p_0$ of the Hecke algebra $\h(G,U)$ do not coincide, i.e. the canonical surjection
\begin{align*}
C^*(L^1(G, U)) \longrightarrow p_0 C^*({G})p_0
\end{align*}
is not injective.
\end{cor}

A result of similar flavour was established in \cite{valette}; explicitly, it was shown right after \cite[Corollary 16]{valette} that the canonical $*$-homomorphism from $C^*(L^1(G, K))$ to $C^*(G)$ is not surjective for any non-compact simple Lie group $G$ of rank $\geq 2$ with maximal compact subgroup $K$.

\appendix\label{appendix:T}

\section{Tzanev's property T and relation to the Schlichting completion}

In \cite[Chapitre 1]{tzanev-thesis}, Tzanev defined a property (T) for a Hecke pair $(G, \Gamma)$, and claimed without proof that a Hecke pair $(G, \Gamma)$ has property (T) if and only if the Schlichting completion $\overline{G}$ has property (T) in the ordinary sense for groups. He further claimed that all other characterisations of property (T) translate to Hecke pairs. We have already seen that the last claim is more subtle, since isolation of the trivial representation fails when regarded at the level of the $C^*$-completion of $L^1(G, U)$ for $(G, U)$ as in Theorem~\ref{SLn trivial representation not isolated}.

Tzanev's first claim is true, but in filling out the details of the proof we have found a more involved argument than expected. We include the proof in this appendix, in the hope that the result could be useful elsewhere.

 Recall the definition of property (T) for locally compact groups, \cite{BdHV}. Let $\mathscr{G}$ be a locally compact group and $\pi:\mathscr{G} \to \mathcal{U}({H})$ a unitary representation on a Hilbert space $H$. Given $\varepsilon >0$ and $\mathscr{Q}$ a compact subset of $\mathscr{G}$, a vector $\xi\in H$ with $\Vert \xi\Vert=1$ is
 $(\varepsilon, \sq)$-invariant if $\| \pi(g)\xi - \xi \| < \varepsilon$ for all $g\in \sq$. The representation $\pi$ has almost $\sg$-invariant vectors if it admits an $(\varepsilon, \sq)$-invariant vector for any $\varepsilon>0$ and $\sq$ compact subset of $\sg$.  The group $\sg$ has property (T) of Kazhdan if every unitary representation of $\sg$ having almost $\sg$-invariant vectors has a nontrivial $\sg$-invariant vector.

By analogy with the above, Tzanev introduced a notion of property (T) for a Hecke pair $(G, \Gamma)$. We shall phrase the definition in terms of a topological Hecke pair $(\sg, \sh)$. In fact, the definition makes sense only assuming that $\sh$ is a closed subgroup of $\sg$.

\begin{df}\label{property T for pairs df}
 Let $\mathscr{G}$ be a topological group, $\mathscr{H}$ a closed subgroup and $\pi:\mathscr{G} \to \mathcal{U}({H})$ a unitary representation of $\sg$ on a Hilbert space $H$. Given $\varepsilon >0$ and $\mathscr{Q}$ a compact subset of $\mathscr{G}/\sh$,  a vector $\xi$ in $H$ with  $\Vert \xi\Vert=1$ is $(\varepsilon, \sq)$-invariant if $\xi$ is $\sh$-invariant and $\| \pi(g)\xi - \xi \| < \varepsilon$ for all $[g] \in \sq$.

 The representation $\pi$ has \emph{$\sh$-invariant} \emph{almost $\sg$-invariant vectors} if for every compact subset $\sq \subseteq \sg/\sh$ and every $\varepsilon >0$ there exists a $(\varepsilon, \sq)$-invariant vector.

 The pair $(\sg, \sh)$ has \emph{property (T)} if every unitary representation of $\sg$ having $\sh$-invariant almost $\sg$-invariant vectors has a nontrivial $\sg$-invariant vector.\\
\end{df}

In the case of a discrete Hecke pair $(G, \Gamma)$, the above recovers Tzanev's definition upon replacing $\sq$ with  a finite subset $K/\Gamma$ of $G/\Gamma$ for $K$ a subset of $G$.

\begin{rem}
 The above notion of property (T) is not the same as the relative property (T) for pairs. On one hand, $(\sg, \sg)$ always has property (T) in the sense of Definition \ref{property T for pairs df}, but only has the relative property (T) when $\sg$ has property (T). On the other hand, $(\sg, \{e\})$ always has the relative property (T), but only has property (T) in the sense of Definition \ref{property T for pairs df} when $\sg$ has property (T).
\end{rem}

\begin{rem}\label{rm:alternative-propT} The referee has suggested to us a (natural) equivalent formulation of property (T) for a topological pair $(\sg, \sh)$: the trivial representation $1_{\sg}$ is isolated in  $\mathcal{R}\cup \{1_{\sg}\}$, where $\mathcal{R}$ is any set of equivalence classes of unitary representations of $\sg$ that have nontrivial $\sh$-invariant vectors but do not have nontrivial $\sg$-invariant vectors. Here $\mathcal{R}$ is endowed with its Fell topology, see \cite[Definition F.2.1]{BdHV}. Since all representations in $\mathcal{R}$ have  $\sh$-invariant vectors, the open basic sets for the topology may be characterised by compact subsets $\sq$ of $\sg/\sh$. Then the claimed equivalence can be proved similar to \cite[Proposition 1.2.3]{BdHV}.
\end{rem}

\begin{prop}
\label{property T for normal subgroups}
 Let $\sg$ be a topological group and $ \sn \subseteq \sh \subseteq \sg$  closed subgroups, with $\sn$ a normal subgroup of $\sg$. The pair $(\sg, \sh)$ has property (T) if and only if  $(\sg / \sn, \sh / \sn)$ has property (T).

In particular, for a topological group $\sg$ and a closed normal subgroup $\sn \unlhd \sg$, we have that $(\sg, \sn)$ has property (T) if and only if $\sg / \sn$ has property (T).
\end{prop}

{\bf \emph{Proof:}} $(\Longrightarrow)$ Suppose $(\sg, \sh)$ has property (T). Let $\pi: \sg / \sn \to \mathcal{U}({H})$ be a unitary representation that has $\sh/\sn$-invariant almost $\sg/\sn$-invariant vectors. We want to show that $\pi$ has a nontrivial $\sg /\sn$-invariant vector. Let $\widetilde{\pi}:= \pi \circ q$ be the lifting of $\pi$ to $\sg$ through the quotient map $q :\sg \to \sg / \sn$. Let $\sq \subseteq \sg / \sh$ be a compact set and $\varepsilon > 0$. Under the canonical homeomorphism between $\sg / \sh$ and $({\sg / \sn})/({\sh/\sn})$ we can naturally identify $\sq$ with a compact subset $\sq' \subseteq ({\sg / \sn})/({\sh/\sn})$. Thus, there exists a $\sh/ \sn$-invariant vector $\xi \in {H}$ of norm one such that
\begin{align}
\label{(G,H) has T implies G/H has T equation}
 \|\pi(g\sn) \xi - \xi \| < \varepsilon
\end{align}
 for every $g\sn$ such that $[g\sn] \in \sq' \subseteq ({\sg / \sn})/({\sh/\sn})$. The following diagram of canonical maps commutes:
\begin{align}
\label{commutative diagram G, GN, GH, GN HN}
\xymatrix{\sg \ar[r] \ar[d] & \sg / \sn \ar[d] \\
\sg / \sh \ar[r]^{\cong} & \frac{\sg / \sn}{\sh / \sn}}
\end{align}
Therefore the set of elements $g \in \sg$ such that $[g\sn] \in \sq'$ is the same as the set of elements $g \in \sg$ such that $g\sh \in \sq$. Thus, condition (\ref{(G,H) has T implies G/H has T equation}) simply says that $\| \widetilde{\pi}(g) \xi - \xi \| < \varepsilon$ for all $g \in \sg$ such that $g\sh \in \sq$. It is moreover clear that $\xi$ is $\sh$-invariant for $\widetilde{\pi}$. Hence, $\widetilde{\pi}$ has $\sh$-invariant almost $\sg$-invariant vectors. From property (T) it then follows that $\widetilde{\pi}$ has a true $\sg$-invariant vector, say $\xi_0 \in {H}$. It is clear that $\xi_0$ is then a  $\sg/\sn$-invariant vector for $\pi$.

$(\Longleftarrow)$ Suppose $(\sg / \sn, \sh/ \sn)$ has property (T). Let $\pi: \sg  \to \mathcal{U}({H})$ be a unitary representation that has $\sh$-invariant almost $\sg$-invariant vectors. We want to show that $\pi$ has a true $\sg$-invariant vector. Let us consider the subspace ${H}^{\sn}$ of $\sn$-invariant vectors, i.e. ${H}^{\sn} := \{ \eta \in {H}: \pi(h) \eta = \eta\,, \forall h \in \sn \}$, which is nontrivial because $\pi$ has $\sh$-invariant (hence also $\sn$-invariant) vectors. Consider now the unitary representation $\sigma: \sg / \sn \to \mathcal{U}({H}^{\sn})$ defined by $\sigma([g])\eta := \pi(g)\eta$. Note that this is a representation on ${H}^{\sn}$, i.e. $\pi(g)\eta \in {H}^{\sn}$ whenever $\eta \in {H}^{\sn}$, because normality of $\sn$ in $\sg$ yields that
\begin{align*}
 \pi(h) \pi(g) \eta = \pi(g)\pi(g^{-1} h g) \eta = \pi(g) \eta\,.
\end{align*}
We claim that $\sigma$ has $\sh / \sn$-invariant almost $\sg / \sn$-invariant vectors. Let $\sq \subseteq (\sg / \sn)/(\sh / \sn)$ be a compact set and $\varepsilon > 0$. Under the homeomorphism between $\sg/\sh$ and $(\sg / \sn)/(\sh / \sn)$ we can identify $\sq$ with a compact subset $\sq' \subseteq \sg / \sh$. Thus, there exists an $\sh$-invariant vector $\xi \in {H}$ such that
\begin{align}
\label{G/H has T implies (G,H) has T equation}
 \|\pi(g) \xi - \xi \| < \varepsilon
\end{align}
for all $g \in \sg$ such that $g\sh \in \sq'$. Being $\sh$-invariant (hence $\sn$-invariant) for $\pi$ means that $\xi \in {H}^{\sn}$ and moreover it implies that $\xi$ is $\sh/\sn$-invariant for $\sigma$. Condition (\ref{G/H has T implies (G,H) has T equation}) simply says that $\| \sigma([g])\xi - \xi \| < \varepsilon$ for all $g \in \sg$ such that $g\sh \in \sq'$. By  commutativity of the diagram (\ref{commutative diagram G, GN, GH, GN HN}), it follows that $\| \sigma([g])\xi - \xi \| < \varepsilon$ for all $g \in \sg$ such that $[g\sn] \in \sq$. Thus, $\sigma$ has $\sh / \sn$-invariant almost $\sg / \sn$-invariant vectors and therefore, by property (T) for the  pair $(\sg / \sn, \sh / \sn)$, it must have a nontrivial $\sg / \sn$-invariant vector $\xi_0 \in {H}^{\sn}$. It is clear that $\xi_0$ is then a $\sg$-invariant vector for $\pi$.

The second claim of this proposition, that $(\sg, \sn)$ has property (T) if and only if $\sg/ \sn$ has Property (T), follows directly from the first claim of the proposition by taking $\sh = \sn$. {}

The following result generalises the known fact that Property (T) passes to quotients.

\begin{prop}
\label{(G,K) has T implies (G,H) has T}
 Let $\sg$ be a topological group and $\sk \subseteq \sh \subseteq \sg$ closed subgroups. If $(\sg, \sk)$ has property (T), then $(\sg, \sh)$ has property (T).

In particular, if $\sg$ has property (T), then $(\sg, \sh)$ has property (T) for any closed subgroup $\sh \subseteq \sg$.
\end{prop}

 We shall need the following result, whose routine proof we omit.

\begin{lemma}
 \label{lemma about natural map G/K to G/H}
Let $\sg$ be a topological group and $\sk \subseteq \sh \subseteq \sg$ closed subgroups. The natural map $\psi: \sg / \sk \to \sg / \sh$ defined by $\psi(g\sk) := g\sh$ is continuous.
\end{lemma}

{\bf \emph{Proof of Proposition \ref{(G,K) has T implies (G,H) has T}:}} We assume first that $(\sg, \sk)$ has property (T). Let $\pi: \sg \to \mathcal{U}({H})$ be a unitary representation with $\sh$-invariant almost $\sg$-invariant vectors. We must show that $\pi$ has a true nonzero $\sg$-invariant vector. That will follow immediately from property (T) of $(\sg, \sk)$ if we prove that $\pi$ has $\sk$-invariant almost $\sg$-invariant vectors, which we will now show.

 Let $\sq \subseteq \sg / \sk$ be any compact set and let $\varepsilon > 0$. Then $\psi(\sq)$ is a compact set of $\sg / \sh$ by Lemma \ref{lemma about natural map G/K to G/H}. Therefore, there exists a $\sh$-invariant vector $\xi \in {H}$ of norm one such that $\|\pi(g)\xi - \xi \| < \varepsilon$ for all $[g] \in \psi(\sq)$. Since $\xi$ is $\sh$-invariant, it is also $\sk$-invariant. Moreover, by definition of the map $\psi$, we have that $\|\pi(g)\xi - \xi \| < \varepsilon$ for all $[g] \in \sq$. Thus, $\pi$ has $\sk$-invariant almost $\sg$-invariant vectors.

 The second claim of the proposition follows by taking $\sk := \{e\}$.
 {}

\begin{prop}
\label{G has T iff G with compact has T}
 Let $\sg$ be a topological group and $\sh$ a compact subgroup. Then $\sg$ has property (T) if and only if the pair $(\sg, \sh)$ has property (T).
\end{prop}

{\bf \emph{Proof:}} The left to right  direction follows directly from Proposition \ref{(G,K) has T implies (G,H) has T}.

For the other implication, assume that the pair $(\sg, \sh)$ has property (T). Let $\pi: \sg \to \mathcal{U}({H})$ be a unitary representation that has almost $\sg$-invariant vectors. We must prove that $\pi$ has nontrivial $\sg$-invariant vectors. This will follow immediately from property (T) of $(\sg, \sh)$ if we prove that $\pi$ has $\sh$-invariant almost $\sg$-invariant vectors. However, this can be achieved using the compactness of  $\sh$. Wo leave the details to the reader.
{}\\

The previous result can also be proved using  Remark~\ref{rm:alternative-propT} (we thank the referee for suggesting this argument). Since $\sh$ is compact, integration over $\sh$ with respect to Haar measure gives rise to an idempotent element $P$ in the multiplier algebra of $C^*(\sg)$. Therefore the set of equivalence classes of unitary representations of $\sg$ having $\sh$-invariant vectors is characterised by $\pi(P)\neq 0$, thus this set is open in the set of unitary representations of $\sg$. Hence isolation of $1_{\sg}$ in either one of these sets happens simultaneously.\\

We are finally ready to verify Tzanev's claim that for a discrete Hecke pair $(G, \Gamma)$, property (T) from Definition~\ref{property T for pairs df} is equivalent to property (T) for the topological group $\overline{G}$ in the Schlichting completion  $(\overline{G}, \overline{\Gamma})$. Recall that $(G, \Gamma)$ is \emph{reduced} if $R^\Gamma=\cap_{g\in G}g\Gamma g^{-1}$ is the trivial subgroup $\{e\}$. For an arbitrary Hecke pair  $(G, \Gamma)$, the reduction $(G_r, \Gamma_r):=(G/R^\Gamma, \Gamma/R^\Gamma)$ is a reduced Hecke pair, and there are canonical isomorphisms between $\Hh(G, \Gamma)$, $\Hh(G_r, \Gamma_r)$ and $\Hh(\overline{G}, \overline{\Gamma})$, see \cite{tzanev} and \cite{kal-land-qui}. A different proof of the equivalence of i) and iv) in the next theorem can be found in \cite[Proposition 6.4]{AD}.

\begin{thm}\label{thm:equiv-T-Heckepair}
 Let $(G, \gm)$ be a Hecke pair. The following are equivalent:
\begin{itemize}
 \item[i)] the pair $(G,\gm)$ has property (T);
 \item[ii)] the pair $(G_r, \gm_r)$ has property (T);
 \item[iii)] the pair $(\overline{G}, \overline{\gm})$ has property (T);
 \item[iv)] the group $\overline{G}$ has property (T).
\end{itemize}
Here $G$ and $G_r$ are assumed to have the discrete topology and $\overline{G}$ is assumed to have its locally compact totally disconnected topology.
\end{thm}

{\bf \emph{Proof:}} $i) \Longleftrightarrow ii)$ Follows directly from Proposition \ref{property T for normal subgroups} since $R^\Gamma$ is a normal subgroup of $G$ contained in $\Gamma$.

$iii) \Longleftrightarrow iv)$ Follows directly from Proposition \ref{G has T iff G with compact has T}.

$ii) \Longrightarrow iii)$ Suppose that $(G_r, \gm_r)$ has property (T). Let $\pi: \overline{G} \to \mathcal{U}({H})$ be a unitary representation that has $\overline{\gm}$-invariant almost $\overline{G}$-invariant vectors. We claim that $\pi$ has a true nontrivial $\overline{G}$-invariant vector.

Let $\epsilon > 0$ and $Q \subseteq G_r / \gm_r$ a finite set, say $Q = \{g_1 \gm_r, \dots, g_n\gm_r\}$. Since $\pi$ has $\overline{\gm}$-invariant almost $\overline{G}$-invariant vectors, there is a $\overline{\gm}$-invariant vector $\xi \in {H}$ of norm one such that
\begin{align*}
 \| \pi(g) \xi - \xi \| < \epsilon\,,
\end{align*}
for all $[g] \in \{g_1 \overline{\gm}, \dots, g_n \overline{\gm} \}$. By restriction, the same holds for all $[g] \in \{g_1 \gm_r, \dots, g_n \gm_r \}$. Moreover, being $\overline{\gm}$-invariant, the vector $\xi$ is $\gm_r$-invariant. Hence, we showed that the restriction of $\pi$ to $G_r$ has $\gm_r$-invariant almost $G_r$-invariant vectors. By property (T) for $(G_r, \Gamma_r)$ it follows that there exists a  nontrivial $G_r$-invariant vector. By continuity of $\pi$, this vector must be $\overline{G}$-invariant.

$iii) \Longrightarrow ii)$ Suppose $(\overline{G}, \overline{\gm})$ has property (T). Let $\pi: G_r \to \mathcal{U}({H})$ be a unitary representation that has $\gm_r$-invariant almost $G_r$-invariant vectors. We must show that it has a nontrivial  $G_r$-invariant vector.

Let $\mathscr{V}:= \overline{\pi(G_r) {H}^{\gm_r}}$, where ${H}^{\gm_r}$ is the subspace of $\gm_r$-fixed vectors, which is nontrivial because $\pi$ is assumed to have nontrivial $\gm_r$-invariant vectors. It is clear that $\mathscr{V}$ is an invariant subspace for $\pi$, so the restriction of $\pi$ to this subspace gives rise to a unitary representation $\pi|$ of $G_r$ on $\mathscr{V}$. Moreover, $\pi|$ is clearly generated by the $\gm_r$-fixed vectors. By \cite[Proposition 6.17]{kal-land-qui},  the representation $\pi|$ is continuous with respect to the Hecke topology, and therefore extends uniquely to a representation $\widetilde{\pi|}$ of $\overline{G}$ on $\mathscr{V}$, which is generated by its $\overline{\gm}$-fixed vectors.

Let $\varepsilon >0$ and $Q \subseteq \overline{G} / \overline{\gm}$ be a compact set. Since $\overline{G} / \overline{\gm}$ is discrete, $Q$ is in fact a finite set $Q = \{g_1 \overline{\gm}, \dots, g_n\overline{\gm} \}$, where we can assume without loss of generality that the representatives $g_1, \dots, g_n$ are elements of $G_r$, because of the canonical bijection between $\overline{G} / \overline{\gm}$ and $G_r /\gm_r$. For this chosen $\varepsilon$ and the finite set $\widetilde{Q}:= \{ g_1 \gm_r, \dots, g_n \gm_r \}$, there exists a $\gm_r$-invariant unit vector $\xi \in \mathscr{V}$ such that
\begin{align*}
 \| \pi|(g)\xi - \xi \| < \varepsilon\,,
\end{align*}
for every $[g] \in \widetilde{Q} \subseteq G_r / \gm_r$. By continuity, the vector $\xi$ is $\overline{\gm}$-invariant for $\widetilde{\pi|}$ and we have that $\| \widetilde{\pi|}(g)\xi - \xi \| < \varepsilon$ for every $[g] \in Q$. Hence, $\widetilde{\pi|}$ has $\overline{\gm}$-invariant almost $\overline{G}$-invariant vectors. By property (T) for $(\overline{G}, \overline{\gm})$, the representation $\widetilde{\pi|}$  has a nontrivial $\overline{G}$-invariant vector. Thus, by restricting $\widetilde{\pi|}$ to $G_r$, it follows that $\pi$ has a nontrivial $G_r$-invariant vector. {}

\end{document}